\documentclass[11pt]{amsart}
\usepackage{amssymb}
\usepackage{amsmath}

\usepackage{hyperref}
\usepackage{graphicx,psfrag,epsfig,multirow}

\usepackage[latin1]{inputenc}
\usepackage{enumerate}

\usepackage{color}

\definecolor{orange}{rgb}{1,0.5,0}

\def\Z{{\mathbb Z}}\def\T{{\mathbb T}}\def\R{{\mathbb R}}\def\C{{\mathbb C}}

\def\cC{\mathcal{C}}
\def\one{\mathbf{1}}

\def\o{\omega}\def\e{\eta}

\def\al{\alpha}\def\be{\beta}
\def\be{\beta}
\def\th{\theta}

 {}

\def\a#1{\left|#1\right|}

\def\<{\langle}
\def\>{\rangle}

\theoremstyle{plain}

\def\cC{\mathcal C}
\def\cM{\mathcal M}

\def\a{\alpha}

\def\e{\varepsilon}
\def\pa{\partial}

\def\th{\theta}

\let\newpf\proof \let\proof\relax

\def\cT{\mathcal{T}}

\newcommand{\ba}{\overline{A}}

\newcommand{\cF}{\mathcal{F}}

\newcommand{\cO}{\mathcal{O}}

\def\be{\begin{equation}}
\def\ee{\end{equation}}

\def\ba{{\begin{align}}}
\def\ea{{\end{align}}}

\def\bm{\begin{pmatrix}}
\def\em{\end{pmatrix}}

\def\a{{\alpha}}

\def\bS{\mathbb{S}}
\def\bD{\mathbb{D}}

\def\0{{\mathbf 0}}

\newtheorem{question}{Question}
\newtheorem{thm}{Theorem}[section]

\newtheorem{conjecture}[thm]{Conjecture}

\theoremstyle{remark}

\theoremstyle{definition}

\newenvironment{proof}{ \noindent{\it Proof.}\quad}{\ \hfill $\Box$\vskip .2cm}

\newcommand{\eps}{{\epsilon}}

\newcommand{\om}{{\omega}}

\newcommand{\N}{{\mathbb N}}
\newcommand{\Q}{{\mathbb Q}}

\def\B0{{\bold{0}}}

\def\a{\alpha}

\def\be{\begin{equation}}
\def\ee{\end{equation}}

\def\cC{\mathcal{C}}
\def\cF{\mathcal{F}}

\def\1{{\bf 1}}

\catcode`\@=12

\def\Empty{}
\newcommand\oplabel[1]{
  \def\OpArg{#1} \ifx \OpArg\Empty {} \else
  	\label{#1}
  \fi}
		
\newcommand{\comm}[1]{}
\newcommand{\comment}[1]{}

\begin{document}
\title{Some questions around quasi-periodic dynamics}
\author{Bassam Fayad and  Raphaël Krikorian}
\address{Bassam Fayad, IMJ-PRG CNRS}
\email{bassam@math.jussieu.fr} 
\address{Raphaël Krikorian, Universit\'{e} de Cergy-Pontoise} 
\email{raphael.krikorian@u-cergy.fr}
\thanks{BF and RK are supported by project ANR BEKAM: ANR-15-CE40-0001. RK is supported by a Chaire d'Excellence LABEX MME-DII }
\maketitle

We propose in these notes a list of some old and new questions related to quasi-periodic dynamics. A main aspect of quasi-periodic dynamics is the crucial influence of  arithmetics on the dynamical features, with a strong duality in general between Diophantine and Liouville behavior. We will discuss rigidity and stability in Diophantine dynamics as well as their absence in Liouville ones. Beyond this classical dichotomy between the Diophantine and the Liouville worlds, we discuss some unified approaches and some phenomena that are valid in both worlds. Our focus is mainly on low dimensional dynamics such as circle diffeomorphisms, disc dynamics, quasi-periodic cocycles, or surface flows, as well as finite dimensional Hamiltonian systems. 

In an opposite direction, the study of the dynamical properties of some  diagonal and unipotent actions on the space of lattices can be applied to arithmetics, namely to the theory of Diophantine approximations. We will mention in the last section some problems related to that topic.

The field of quasi-periodic dynamics is very extensive and has a wide range of interactions with other mathematical domains. The list of questions we propose is naturally far from exhaustive and our choice was often motivated by our research involvements. 
 
 \section{Arithmetic  conditions}  A vector $\a\in\R^d$ is {\it non-resonant} if it  has rationally independent coordinates:  for all $(k_{1}\ldots, k_{d})\in\Z^{d}$, the identity $\sum_{i=1}^d k_{i}\a_{i}=0$ implies $k_{i}=0$ for $i=1,\ldots, d$; otherwise, it is called {\it resonant}. 

For $\gamma,\sigma>0$, we define the set $DC_{d}(\gamma,\sigma)\subset\R^d$ of {\it diophantine} vectors  with {\it exponent} $\sigma$ and {\it constant} $\gamma$  as the set of  $\a=(\a_{1},\ldots,\a_{d})\in\R^d$ such that 
\be \forall (k_{1},\ldots, k_{d})\in\Z^{d},\ |\sum_{i=1}^d k_{i}\a_{i}|\geq \frac{\gamma}{(\sum_{i=1}^d|k_{i}|)^\sigma};
\ee
we then set  $DC_{d}(\sigma)=\bigcup_{\gamma>0}DC(\gamma,\sigma)$, $DC_{d}=\bigcup_{\sigma>0}DC(\sigma)$.   For each fixed $\sigma>d$ and $\gamma$ small enough the set $DC(\gamma,\sigma)$ has positive Lebesgue measure in the unit ball of $\R^d$ and the Lebesgue measure of its complement goes to zero as $\gamma $ goes to zero. Thus the sets $DC(\sigma)$, $\sigma>d$ and  $DC_{d}$ have full Lebesgue measure in $\R^d$. The set $DC_{d}$ is the set of   {\it Diophantine}  vectors of $\R^d$  while  its complement in the set of non-resonant vectors  is called the set of  {\it Liouville}  vectors.

 For 	a translation vectors of $\T^d$ defined as  $\a+\Z^d$, $\a \in\R^d$, we say that it is resonant, Diophantine or Liouville, if the $\R^{d+1}$ vector $(1,\a)$ is  resonant, Diophantine or Liouville respectively.

\section{Diffeomorphisms of the circle and the torus}

For $k\in\N\cup\{\infty,\omega\}$ we define  ${\rm Diff}_{0}^k(\T^d)$   as the set of  orientation preserving homeomorphisms  of $\T^d$ of  class $C^k$ together with their inverse. To any $f\in {\rm Diff}_{0}^0(\T^d)$ one can associate its  {\it rotation set}   $\rho(f):=\{\int_{\T}(\bar f-id)d\mu,\mu\in\cM(f)\}  \mod \Z^d$ where $\bar f:\R^d\to\R^d$ is a lift of $f$ and $\cM(f)$ is the set of all $f$-invariant probability measures on $\T^d$. 
Let  $T_{\a}:\T^d\to\T^d$ be the translation $x\mapsto x+\a$, $\cF^k_{\a}(\T^d)=\{f\in{\rm Diff}_{0}^k(\T^d),\ \rho(f)=\{\a\}\}$, $\cO_{\a}^k(\T^d)=\{h\circ T_{\a}\circ h^{-1}, h\in{\rm Diff}_{0}^k(\T^d)\}$. We say that $f\in {\rm Diff}_{0}^\infty(\T^d)$ is {\it almost reducible} if there exists a sequence $(h_{n})_{n\in\N}\in ({\rm Diff}_{0}^\infty(\T^d))^\N$ such that $h_{n}\circ f\circ h_{n}^{-1}$ converges in the $C^\infty$-topology  to $T_{\a}$.
When $d=1$,  $\rho(f)$ is reduced to a single element  and we denote by $\rho(f)$ this element. By Denjoy Theorem, any $f\in{\rm Diff}^k_{0}(\T)$ with $k\geq 2$,  is conjugated by an orientation preserving  homeomorphism to $T_{\a}$. If furthermore $\a$ is Diophantine and $k=\infty$ then by Herman-Yoccoz theorem \cite{H-ihes}, \cite{Y-smoothdiffeo} this conjugacy is smooth which amounts to $\cF_{\a}^\infty(\T)=\cO_{\a}^\infty(\T)$. It is of course natural to try to extend this result to the higher dimensional situation where $f$ is an orientation preserving diffeomorphism of the $d$-dimensional torus $\T^d$. Unfortunately,  no Denjoy theorem is available in this situation and the only reasonable question to ask for  is the following 
\begin{question}\label{diffS1:red} Let $f:\T^d\to\T^d$ be a smooth diffeomorphism of the torus $\T^d=\R^d/\Z^d$ which is topologically conjugate to a translation $T_{\a}:\T^d\to\T^d$, $x\mapsto x+\a$ with $\a$ Diophantine. Is  the conjugacy smooth?
\end{question}
Notice that when $d=2$, even if $\a$ is Diophantine,  $\cF_{\a}^\infty(\T^2)$ is not equal to $\cO_{\a}^\infty(\T^2)$ or $\overline{\cO_{\a}^\infty(\T^d)}$ as is shown by taking projectivization of cocyles in $SW^\infty(\T,SL(2,\R))$: such cocycles have a uniquely defined rotation number, that can be chosen Diophantine, and at the same time  can have positive Lyapunov exponents (which prevents the projective action to be conjugated to a translation) (cf. \cite{Herman83}). Analogously, by taking projectivization of cocycles in $SW^\omega(\T,SL(2,\R))$ and using Avila's theory characterizing   sub-critical/critical cocycles and the Almost Reducibility Conjecture  (see Section \ref{sec:5}) one can  show that there exist elements of $\overline{\cO_{\a}^\infty(\T^d)}$ which are not $C^\infty$- almost reducible and, even if $\a\in\T^d$ is Diophantine, that the set $\cO_{\a}^\infty(\T^d)$ is not  closed.

In a similar vein
\begin{question}\label{diffS1:ared} Let $f:\T^d\to\T^d$ be a smooth diffeomorphism  which is topologically conjugate to the translation with $\a$ non-resonant. Is it $C^\infty$-accumulated by elements of $\cO_{\a}^\infty(\T^d)$? Is  it  $C^\infty$-almost reducible? 
\end{question}
When $d=1$ the  first and the second  part of the  preceding question have a positive answer. Yoccoz proved \cite{Y-centralizers} that $\cF_{\a}^\infty(\T)=\overline{\cO_{\a}^\infty}(\T)$ and it  is proved in \cite{AK4} that any smooth orientation preserving diffeomorphism of the circle is $C^\infty$-almost reducible. The proof of this result   uses renormalization techniques which at the present time doesn't seem to extend to the higher dimensional case.  Still the situation in the semi-local case might be more accessible.

\begin{question} \label{Q:semilocs1}Same questions as in Question \ref{diffS1:red} - \ref{diffS1:ared} in the semi-local case that is for $f$ in some neighborhood of the set of rotations, independent of $\a$.
\end{question}
If   one assumes $\a$ to be Diophantine and $f$ to   be in a neighborhood of $T_{\a}$ that {\it depends} on $\a$ the answer to Question \ref{diffS1:red} is positive; this can be proved by standard KAM techniques.

\section{Pseudo-rotations of the disc}
A $C^k$ ($k\in\N\cup\{\infty,\omega\}$) {\it pseudo rotation} of the disk $\bD=\{(x,y)\in\R^2, x^2+y^2\leq 1\}$ is a $C^k$ orientation and {\it area preserving }diffeomorphism of the disk $\bD$ that fixes the origin, leaves invariant the boundary $\pa\bD$ of the disk and with {\it no other periodic point than the origin}. 
 Like in the case of circle diffeomorphisms one can define for such pseudo-rotation a unique rotation number around the origin which is invariant by conjugation  (see for example \cite[Corollary 2.6]{franks.etds} or \cite[Theorem 3.3]{franks}). 
Anosov and Katok \cite{AK} constructed  in 1970, {\it via}  approximation by periodic dynamics, {\it ergodic} (for the area measure)  and infinitely differentiable pseudo-rotations of the disk, providing thus the first examples of pseudo-rotations  which are not topologically conjugate to rigid rotations.  By a theorem of Franks and Handel \cite{franks-handel} a {\it transitive} area and orientation preserving diffeomorphism of the disk fixing the origin and leaving invariant the boundary of the disk must be a pseudo-rotation. 
\subsection{Birkhoff rigidity conjecture}
A famous question on pseudo-rotations attributed to Birkhoff is the following. 
\begin{question} Is a real analytic pseudo-rotation of angle $\a$ analytically conjugated to the rotation $R_{\a}$ of angle $\a$ on the disc? 
\end{question}
Addressing this question should involve the artihmetics of $\a$. On one hand,
R\"{u}ssmann \cite{russmann}  proved  the following alternative  for a Diophantine (in fact of Brjuno type is sufficient) elliptic fixed point of a real analytic  area preserving surface diffeomorphism $f$: either the point is surrounded by a positive measure set of invariant circles with different Diophantine frequencies, or the map $f$ is locally conjugate to a rotation in the neighborhood of the fixed point.  
On the other hand, when the real analytic category is relaxed to infinite differentiability, Anosov-Katok construction provides many counter-examples to the preceding question (for Liouville $\a$'s). 
We can thus divide the preceding question into two questions 
\begin{question}\label{q:anpr} Can one construct Anosov-Katok examples ({\it viz.} ergodic pseudo-rotations) in the real analytic category? If possible, can one impose the rotation number to be any non-Brjuno number?
\end{question}\begin{question}[Reducibility]\label{Qredpd}Is it true that every $C^k$, $k=\infty,\omega$, pseudo-rotation of the disk with diophantine rotation number $\a$ $C^k$-conjugated to a rigid rotation by angle $\a$?
\end{question}
Notice that in the smooth category the answer to Question \ref{q:anpr}  is positive: for the first part this is the existence of Anosov-Katok ergodic, even weak mixing, pseudo-rotations and for the second part one can prove  that for any Liouville number $\a$, there exists weak mixing pseudo-rotations  as well as examples that are isomorphic to the rotation of frequency $\a$ on the circle \cite{FS,FSW}. Together with Herman's last geometric theorem, this gives in the $C^\infty$-case a complete dichotomy between Diophantine and Liouville behavior.

Let $\cF_{\a}^\infty$ be the set of $C^\infty$ pseudo-rotations with rotation number $\a$ and $\cO_{\a}^\infty$ be the set of $h\circ R_{\a}\circ h^{-1}$ where $h$ is a $C^\infty$ area and orientation preserving map of the disk fixing 0 and leaving invariant the boundary of the disk. A weaker question in the smooth case is :
\begin{question}\label{Qredpd'}For $\a$ diophantine is $\cO_{\a}^\infty$ closed for the $C^\infty$-topology?
\end{question}
In fact, a more general question than Question \ref{Qredpd} is the following:
\begin{question}[Almost reducibility]\label{Qaredpd}Is any $C^k$-pseudo-rotation $k=\infty,\omega$, $f$ of the disk with irrational  rotation number $\a$ almost reducible: there exists a sequence of area preserving smooth map $h_{n}$ such that $h_{n}\circ f\circ h_{n}^{-1}$ converges in the $C^k$ topology to $R_{\a}$ (in the analytic case this convergence should occur on a fixed complex neighborhood of the disk)?
\end{question}
Question \ref{Qredpd} has a positive answer in the local case (R\"ussmann for $k=\omega$, Herman, \cite{FKlast} for $k=\infty$) that is when $f$ is in some   $C^k$-neighborhood of $R_{\a}$ (the size of this neighborhood {\it depending} on the  {\it arithmetics} of $\a$). Thus, a positive  answer to   Question \ref{Qaredpd} would imply a positive answer to Question  \ref{Qredpd}. When $k=\infty$,  Question  \ref{Qaredpd} (hence Question  \ref{Qredpd})  has a positive answer in the semi-local case \cite{AK4} that is with the extra assumption that for some $k$ and $\e$  {\it independent } of $\a$, the $C^k$-norm of $Df-id$ is less than $\e$. In this situation one also has 
$\cF_{\a}^\infty\cap W\subset \overline{\cO_{\a}}^\infty$, where $W$ is a neighborhood  for the $C^\infty$-topology of the set of rigid rotations. The proof of the result of \cite{AK4} is based on renormalization techniques and on the fact (proved in \cite{AFLXZ}) that if one has a control on the $C^1$-norm of a pseudo-rotation $f$, the displacement $\max_{\bD}\|f-id\|$ polynomially compares with the rotation number of $f$. Such a control is in general not true for diffeomorphisms of the circle. It is thus natural to ask:

\begin{question}
Describe the set of smooth diffeomorphisms of the circle that are obtained as the restriction on $\bD$ of the dynamics of pseudo-rotations?
\end{question}
\subsection{Rigidity times, mixing and entropy} A diffeomorphism of class $C^k$, $k\in \N \cup \{\infty\}$, is said to admit  {\it $C^k$ rigidity times} (or for short is $C^k$-rigid) if there exists a sequence $q_n$ such that $f^{q_n}$ converges to the Identity map in the $C^k$ topology. 
If we just know that the latter holds in a fixed neighborhood of some point $p$, we say that $f$ is {\it $C^k$ locally rigid} at $p$. All the smooth examples on the disc or the sphere obtained by the Anosov-Katok method are {\it $C^\infty$-rigid} by construction.  
Obviously, rigidity or local rigidity precludes mixing. Hence, the following natural question was raised in \cite{FaKa} in connection with the smooth realization problem and the Anosov-Katok construction method. 

\begin{question} \label{q.FK} Is it true that  a smooth area preserving diffeomorphism of the disc  with zero metric entropy  is not mixing?
\end{question}

In the case of zero topological entropy, and in light of Franks and Handel result, the question becomes
\begin{question}  Is it true that a smooth pseudo-rotation is not mixing?
\end{question}
Bramham \cite{bramham} proved that this is true  if the rotation number is sufficiently Liouville; indeed he proves in that case the existence of $C^0$-rigidity times. It was shown in \cite{AFLXZ} that real analytic pseudo-rotations (with no restriction on the rotation number) are never topologically mixing. By a combination of KAM results and control of recurrence for pseudo-rotations with Liouville rotation numbers, it is actually shown that real analytic pseudo-rotations are  $C^\infty$ locally rigid near their center.

Note that the following is not known, except in $C^1$ regularity where a positive answer is given by \cite{jairo}.
\begin{question}  Does there exist a smooth area preserving disc diffeomorphism that has zero metric entropy and positive topological entropy? 
\end{question}

The following question was raised by Bramham in \cite{bramham}. 
\begin{question} Does every $C^k$ pseudo-rotation $f$ admit $C^0$ rigidity times? The question can be asked for any $k \geq 1$, $k=\infty$ or $k = \o$. 
\end{question}

In the case $k=\o$ or  $\rho(f)$ Diophantine and $k=\infty$, the latter question becomes an intermediate question relative to the Birkhoff-Herman problem on the conjugability of $f$ to the rigid disc rotation of angle $\rho(f)$.
In \cite{AFLXZ} it was shown that for every irrational $\a$, if an analytic  pseudo-rotation of angle $\a$ is sufficiently close to $R_\a$ then it admits $C^\infty$-rigidity times.

\begin{question} Given a fixed analyticity strip, does there exist  $\epsilon>0$ such that if a  real analytic pseudo-rotation  is $\epsilon$ close to the rotation on the given analyticity strip, then it is rigid? 
\end{question}

An {\it a priori} control on the growth of $\|Df^{m}\|$ for a pseudo-rotation is sufficient to deduce the existence of rigidity times  for larger classes of rotation numbers. If for example a polynomial bound holds on the growth of $\|Df^{m}\|$ for a smooth pseudo-rotation, then the existence of  $C^\infty$ rigidity times  would follow for any Liouville rotation number (see \cite{AFLXZ}). In the case of a circle diffeomorphism $f$  a gap in the growth of these norms  is known to hold between exponential growth in the case $f$ has a hyperbolic periodic point or a growth bounded by $O(m^2)$ if not \cite{polterovich}. Does a similar dichotomy hold for area preserving disc diffeomorphisms?

\begin{question} \label{q.liouville} Is there any polynomial bound on the growth of the derivatives of a pseudo-rotation?  Does  every $C^\infty$ pseudo-rotation with Liouville rotation number admit $C^0$ (or even $C^\infty$) rigidity times?
\end{question}

 With Herman's last geometric theorem, a positive answer to the second part of Question \ref{q.liouville} would imply that smooth pseudo-rotations, and therefore  area preserving smooth diffeomorphisms of the disc with zero topological entropy  are never topologically mixing.

In the proof of absence of mixing of an analytic pseudo-rotation, \cite{AFLXZ} uses an {\it a priori} bound on the growth of the derivatives of the iterates of a pseudo-rotation that is obtained {\it via} an 
 effective finite information version of the Katok closing lemma for an area preserving surface diffeomorphism $f$. This effective result  provides a positive  gap in the possible growth of  the derivatives of $f$ between exponential and sub-exponential. 

 In \cite{FZ}, an explicit
finite information condition is obtained for area preserving $C^2$ surface diffeomorphisms,  
that guarantees positive topological entropy.
\begin{question} Find a finite information condition on the complexity growth of an area preserving $C^2$ surface diffeomorphism that insures positive metric entropy. 
\end{question}

Finally, inspired by R\"ussmann and Herman's last geometric theorem on one hand,  and  the Liouville  pseudo-rotations rigidity on the other, we ask the following 
\begin{question} Can a smooth area preserving diffeomorphism of a surface that has an irrational elliptic fixed point be topologically mixing? Can it have an orbit that converges to the fixed point?
\end{question}

\section{Hamiltonian systems}

A $C^2$ function $H:(\R^{2d},0)\to\R$ such that $DH(0)=0$ defines on a neighborhood of 0 a hamiltonian vector field $X_{H}(x,y)=(\pa_{y}H(x,y),-\pa_{x}H(x,y))$  and its flow $\phi^t_{H}$ is a flow of symplectic diffeomorphisms preserving the origin. We shall assume that $0 \in \R^{2n}$ is  an elliptic equilibrium point with $H$ of the following form 
\begin{equation} \label{Hamintro} H(x,y)= \sum_{j=1}^d \omega_j(x_j^2+y_j^2)/2 + O_3(x,y), \end{equation}
where the frequency vector $\omega$ is {\it non-resonant}.

Alternatively we may take $H$ a $C^2$ function
defined on $\T^d \times \R^d$ and consider its Hamiltonian flow $X_{H}(\th,r)=(\pa_{r}H(\th,r),-\pa_{\th}H(\th,r))$. If
\begin{equation} H(\th,r)= \langle \o_0,r\rangle+\cO(r^2) \label{HH} \end{equation} 
then the torus $\T^d\times \{0\}$ is invariant under the Hamiltonian flow
and the induced dynamics on this torus  is the translation
$\phi_{H}^t:\th\mapsto \th+t\o_0.$
Moreover this torus is Lagrangian with respect to the 
canonical symplectic form $d\th\wedge dr$ on $\T^d\times \R^d$. When $\omega$ is Diophantine we say that this torus is a KAM torus.

The stability of an equilibrium or of an invariant quasi-periodic torus by a Hamiltonian flow can be studied from three points of view. 
The usual topological or Lyapunov stability, the stability in a measure theoretic or probabilistic sense which can be addressed by KAM theory  (Kolmogorov, Arnold, Moser), or the effective stability in which one is interested  in quantitative stability in time.

\subsection{Topological stability}

Arnold conjectured that apart from two cases, the case of a sign-definite quadratic part, and generically for $d=2$, an elliptic equilibrium point is generically unstable. 

\begin{conjecture}[Arnold] An elliptic equilibrium point of a generic analytic Hamiltonian system is Lyapounov unstable, provided $n \geq 3$ and the quadratic part of the Hamiltonian function at the equilibrium point is not sign-definite.
\end{conjecture}  

Despite a rich literature and a wealth of results in the $C^\infty$ smoothness (to give a list of contributions would exceed the scope of this presentation), this conjecture is wide open in the real analytic category,  to such an extent that under our standing assumptions (real-analyticity of the Hamiltonian and a non-resonance condition on the frequency vector) not a single example of instability is known.

\begin{question} Give examples of an analytic Hamiltonian that have a non-resonant elliptic equilibrium  (or a non-resonant Lagrangian quasi-periodic torus) that is Lyapunov unstable. 
\end{question}

\begin{question} Give examples of an analytic Hamiltonian that have a non-resonant elliptic equilibrium   (or a non-resonant Lagrangian quasi-periodic torus) that attracts an orbit (distinct from the equilibrium or the torus itself). 
\end{question}
 In \cite{FMS} an example is given of a Gevrey regular Hamiltonien on $\R^6$ that has a non-resonant fixed point at the origin and that has an orbit distinct from the origin that converges to it in the future. 
In \cite{KS12,KG14}, Arnold diffusion methods are used to yield in particular orbits that have $\al$-limit or
$\om$-limit sets that are non-resonant invariant
Lagrangian tori instead of a single non-resonant fixed point.

Following Perez-Marco we ask:
\begin{question} Is it true that a smooth Hamiltonian flow with a non-resonant elliptic equilibrium isolated from periodic points has a hedgehog (a totally invariant compact connected set containing the origin)?
\end{question}

Regarding the additional stability features of elliptic fixed points in the case of two degrees of freedom, we ask the following
\begin{question} Is the iso-energetic twist condition the optimal condition for Lyapunov stability of an irrational elliptic equilibrium in two degrees of freedom?
\end{question}
A smooth example of an irrational equilibrium was constructed by F. Trujillo that satisfies the Kolmogorov non degeneracy condition in $d=2$ degrees of freedom and that has diffusing orbits in some special energy levels.

\subsection{Beyond the classical KAM theory.} 

An equilibrium (or an invariant torus) of a Hamiltonian system is said to be KAM stable if it is accumulated by a positive measure of invariant KAM tori, and if the set of these tori has density one in the neighborhood of the equilibrium (or the invariant torus).

\subsubsection{Weak transversality conditions.}  
 In classical KAM theory, an elliptic fixed point is shown to be KAM-stable under the hypothesis that the frequency vector at the fixed point is non-resonant (or just sufficiently non-resonant) and that the Hamiltonian is sufficiently smooth and satisfies  a generic  non degeneracy condition of its Hessian matrix at the fixed point. Further development of the theory allowed to relax the non degeneracy condition. In \cite{efk_point} KAM-stability was established for non-resonant elliptic fixed points under the (most general) R\"ussmann transversality condition on the Birkhoff normal form of the Hamiltonian. Similar results were obtained for Diophantine invariant tori in \cite{EFK}.

\subsubsection{Absence of transversality conditions.}

\begin{conjecture} \label{herman}[Herman] Prove that an  elliptic equilibrium with a diophantine frequency or a KAM torus  of an analytic Hamiltonian  is accumulated by a set of positive measure of KAM tori.
\end{conjecture}

Clearly, one can of course ask whether KAM stability also holds.

Conjecture \ref{herman} was was made by M. Herman in his ICM98 lecture (in the context of symplectomorphisms). The conjecture is known to be true in two degrees of freedom \cite{russmann}, but remains open in general. It is shown in \cite{EFK}  that an analytic invariant torus $\cT_0$ with Diophantine frequency $\o_0$ is never isolated due to the following alternative.   If the  Birkhoff normal form of the Hamiltonian  at $\cT_0$ satisfies a R\"ussmann transversality condition, the  torus $\cT_0$ is accumulated by KAM tori of positive total measure. If the Birkhoff normal form is degenerate,  there exists a subvariety of dimension at least $d+1$ that is foliated by analytic invariant tori with frequency $\o_0$. 

For Liouville frequencies, one does not expect the conjecture to hold.
\begin{question} Give an example of an analytic Hamiltonian that has a non-resonant (Liouville) elliptic equilibrium that is not is accumulated by a set of positive measure of KAM tori. 
\end{question}

In the $C^\infty$ category (or Gevrey), counter-examples to stability with positive probability can be obtained: in 2 or more degrees of freedom for Liouville frequencies; and in 3 or more degrees of freedom for any frequency vector (\cite{EFK} for $d\geq 4$ and \cite{FS} for $d\geq 3$).  In the remaining case of Diophantine equilibrium with $d=2$, Herman proved stability with positive probability without any twist condition (see \cite{FKlast}).

\subsection{Effective stability}

Combining KAM theory, Nekhoroshev theory and estimates of Normal Birkhoff forms, it was proven in \cite{BFNKAM} that generically, both in a topological and measure-theoretical sense, an invariant Lagrangian Diophantine torus of a Hamiltonian system is doubly exponentially stable in the sense that nearby solutions remain close to the torus for an interval of time which is doubly exponentially large with respect to the inverse of the distance to the torus. It is proven there also  that for an arbitrary small perturbation of a generic integrable Hamiltonian system, there is a set of almost full positive Lebesgue measure of KAM tori which are doubly exponentially stable. These results hold true for real-analytic but more generally for Gevrey smooth systems.  Similar results for elliptic equilibria are obtained in \cite{BFNpoint}.

\begin{question} Give examples of analytic or Geverey differentiable Hamiltonians that have a Diophantine elliptic equilibrium with positive definite twist, that is not more than doubly-exponentially stable in time. Show that this is generic.
\end{question}

\begin{question} Give an example of an analytic Hamiltonian that has a non-resonant elliptic equilibrium with positive definite twist that is not more than exponentially stable in time. 
\end{question}

\begin{question} Give an example of an analytic Hamiltonian that has a  Diophantine elliptic equilibrium that is not more than exponentially stable in time. 
\end{question}

\subsection{On invariant tori of convex Hamiltonians}

\subsubsection{The ''last invariant curve'' of annulus  twist maps.} A classic topic in Hamiltonian systems is that of the regularity of the invariant curves of annulus  twist maps.
 A celebrated result of Birkhoff states that such curves (if they are not homotopic to a point) must be Lipschitz. Numerical evidence seems to indicate that invariant curves are always at least $C^1$. After Mather and Arnaud we ask the following. 
 
\begin{question} Give an example  of a $C^r$, $r\in[2,\infty) \cup \{\omega\},$ annulus twist map that has an invariant $C^0$ but not $C^1$ curve with minimal restricted dynamics. 
\end{question}
In \cite{AF}, a $C^1$ example is constructed, and \cite{arnaud} gives a $C^1$ example with an invariant $C^0$ but not $C^1$ curve having Denjoy type restricted dynamics. 

Due to a result proved by Herman the problem can be reduced to finding  a minimal circle homeomorphism $f$ such that $f+f^{-1}$ is $C^r$ but $f$ is only $C^0$. 

\begin{question} Give an example  of a  $C^r$, $r\in[2,\infty) \cup \{\omega\},$ annulus twist map that has an invariant $C^r$ curve that is not accumulated by other invariant curves. 
\end{question}

\subsubsection{On the destruction of all tori}

Given the Hamiltonian $H=\frac{1}{2} \sum r_i^2$ on $\T^d \times \R^d$. 

\begin{question} What is the maximum of $r$ for which it is possible to perturb $H$ so that the perturbed flow has no invariant Lagrangian torus that is the graph of a  $C^1$ function. 
\end{question}

By Herman, $r \geq d+2-\eps$, $\forall \eps>0$. We also know that $r\leq 2d$ (see \cite{poschel}). In \cite{cheng}, given any frequency $\omega$, a $C^{2d-\eps}$ perturbation of $H$ is given that has no invariant  Lagrangian torus with as unique rotation frequency vector $\omega$.

\subsection{Birkhoff Normal Forms}
Let $H:(\R^{2d},0)\to\R$ be a real analytic  hamiltonian function  admitting 0 as an elliptic non-resonant fixed point. One can always {\it formally} conjugate $H$ to an {\it integrable} hamiltonian: there exist a   {\it formal} (exact) symplectic  germ of diffeomorphism $g$ tangent to the identity and a  {\it formal} series  $N\in \R[[r_{1},\ldots,r_{d}]]$ such that $g_{*}X_{H}=X_{B}$ where $B(x,y)=N(x_{1}^2+y_{1}^2,\ldots,x_{d}^2+y_{d}^2)$. This $B$ is unique and is called the {\it Birkhoff Normal Form} (BNF). This formal object is an invariant of $C^k$-conjugations ($k=\infty,\omega$). Birkhoff Normal Forms can be defined for $C^k$ ($k=\infty,\omega$) symplectic diffeomorphisms admitting an invariant elliptic fixed point or even (in the case of symplectic diffeomorphisms or hamiltonian flows) in a neighborhood of an invariant KAM torus (the frequency must be then diophantine). Siegel \cite{Si54} proved that in general the conjugating transformation could not be convergent and Eliasson asked whether the Birkhoff  Normal Form itself could be convergent.   In the real analytic setting  Perez-Marco \cite{PM} proved that for any given non-resonant quadratic part one has the following dichotomy: either the BNF always converges or it generically diverges. Gong \cite{gong} provided an example of divergent BNF with Liouville frequencies. In \cite{K} it is proved that  the BNF of   a real analytic symplectic diffeomorphism admitting a diophantine elliptic fixed point (with torsion)  is generally divergent.

\begin{question}Let $H$ be a real analytic  Hamiltonian  admitting the origin as a diophantine elliptic fixed point and assume that its Birkhoff Normal Form defines a real analytic function. Is $H$ real analytically conjugated to its Birkhoff Normal Form on a neighborhood of the origin?
\end{question}

\section{Dynamics of quasi-periodic cocycles}\label{sec:5}
Let $G$ be a Lie group (possibly infinite dimensional). A quasi-periodic cocycle  of class $C^k$, $k\in\N\cup\{\infty,\omega\}$ is a map $(\a,A):\T^d\times G\to\T^d\times G$ of the form $(\a,A):(x,y)\mapsto (x+\a,A(x)y)$ where $\a\in\T^d$ (we assume $\a$ to be non-resonant) and $A:\T^d\to G$ is of class $C^k$. We denote the set of  such cocycles $(\a,A)$ by $SW_{}^k(\T^d,G)$ (or $SW_{\a}^k(\T^d,G)$).  The iterates $(\a,A)^n$ of $(\a,A)$ are of the form $(n\a,A^{(n)})$ where  (for $n\geq 1$)   $A^{(n)}$ is the fibered product $A^{(n)}(\cdot)=A(\cdot+(n-1)\a)\cdots A(\cdot+\a)A(\cdot)$.  Two cocycles $(\a,A_{1})$ and $(\a,A_{2})$ are said to be $C^l$-{\it conjugated} if there exists a map $B:\T^d\to G$ (or $B:\R^d/N\Z^d\to G$ for some $N\in\N^*$) of class $C^l$ such that $(\a,A_{2})=(0,B)\circ (\a,A_{1})\circ (0,B)^{-1}$ or equivalently $A_{2}=B(\cdot+\a)A_{1}B(\cdot)^{-1}$. The cocycle $(\a,A)$ is said to be {\it reducible} if it is conjugated to a constant cocycle and, when $H$ is a subgroup of $G$, {\it $H$-reducible} if it is conjugated to an $H$-valued (not necessarily constant) cocycle.  We say that the cocycle is {\it linear} when the group $G$ is a group of matrices. 

\subsection{The case $G=SL(2,\R)$}
 Quasi-periodic $SL(2,\R)$-valued cocycles play an important role in the theory of {\it quasi-periodic Schr\"odinger operators on $\Z$} of the form $H_{x}:l^2(\Z)\to l^2(\Z)$, $H_x: (u_{n})_{n\in\Z}\mapsto (u_{n+1}+u_{n-1}+V(x+n\a)u_{n})_{n\in\Z}$; indeed, the (generalized) eigenvalue equation $H_{x}u=Eu$ 
leads naturally to studying the dynamics of a family of  $SL(2,\R)$-valued  quasi-periodic cocycles depending on $E$, the so-called {\it Schr\"odinger cocycles}.
Many spectral objects or quantities -- such as, resolvent sets (complement  of the spectrum), spectral measures, density of states, speed of decay of Green functions... -- of the family of operators $H_{x}$, $x\in\T^d$, can be related to dynamical notions or invariants for the associated  family of 
Schr\"odinger cocycles -- namely (in that order), uniform hyperbolicity, $m$-functions, fibered rotation number, Lyapunov exponents... We refer to \cite{El-ICM},  \cite{You-ICM} for more details on this topic. 

There are two important quantities associated to $SL(2\,\R)$-valued quasi-periodic cocycles which are invariant by conjugation\footnote{for the rotation number one has to assume the conjugating map to be homotopic to the identity}: the Lyapunov exponent $L(\a,A)$ 
 which measure the exponential speed of growth of the iterates of the cocycle $(\a,A)$ and the fibered rotation number $\rho(\a,A)$ which measures the average speed of rotation of  non-zero vectors in the plane under iteration of the cocycle.
   It is of course tempting to try and classify $SL(2,\R)$-cocycles according to these two invariants.

 The case of real analytic cocycles with one frequency is particularly well understood.  In that situation, following A. Avila \cite{A-global},  one can associate to any cocycle $(\a,A)\in SW_{}^\omega(\T,SL(2,\R))$ a natural family $(\a,A_{\e})\in SW_{}^\omega(\T,SL(2,\C))$ ($\e$ in some neighborhood of 0)  with  $A_{\e}(\cdot)=A(\cdot+\e\sqrt{-1})$. The function $\e\mapsto L(\a,A_{\e})$  plays a very important role in the theory; Avila proved that it is an even  convex  continuous piecewise affine map with {\it quantized} slopes in $2\pi\Z$ (this is the phenomenon of  ``quantization of acceleration'') and that the complex cocycle $(\a,A_{\e})$ is {\it uniformly hyperbolic} if and only  $\e$ is not a break point of  $\e\mapsto L(\a,A_{\e})$.   This analysis leads to the notions of  {\it critical}, {\it supercritical} and  {\it subcritical} cocycles,  where this last term refers to the fact that the function $\e\mapsto L(\a,A_{\e})$ is zero on an neighborhood of $\e=0$. A cocycle $(\a,A)\in SW^\omega_{}(\T,SL(2,\R))$ (homotopic to the identity) can thus have four distinct possible behaviors if one adds to the three preceding ones uniform hyperbolicity. Moreover, the quantization of acceleration allows to {\it predict} the possible transitions between these four regimes and to draw consequences on the spectrum of  Schr\"odinger  operators such as for example the possibility of co-existence of absolutely continuous or pure point spectrum for some type of potentials (cf. \cite{A-global} and for other examples \cite{BK}).  The most striking {\it global} result on the dynamics of these cocycles is certainly the ``Almost reducibility conjecture'' proved by Avila  \cite{A-arac}, \cite{A-kles} which asserts that any {\it subcritical} cocycle in $SW_{}^\omega(\T,SL(2,\R))$ is {\it almost-reducible} (in the analytic category, on a fixed complex neighborhood of the real axis). By \cite{Hou-You}, \cite{You-Zhou}
in the real analytic  semi-local situation (viz. when $A$ is close to a constant, this closeness being independent of $\a$) a cocycle $(\a,A)$ is either uniformly hyperbolic or subcritical. 

   In the $C^\infty$ category, or for many-frequencies systems, our understanding of the dynamics of cocycles is much less complete.  There are important reducibility or almost-reducibility results (\cite{DS}, \cite{El-CMP92}, \cite{K-ast}, \cite{K-annens}, \cite{K-annals}, \cite{AK-annals},\cite{Puig1}, \cite{Puig2}, \cite{rb6}, \cite{AFK}, \cite{Hou-You}, \cite{You-Zhou}, \cite{AK-inventiones}...) but they often involve diophantine conditions and/or  are of perturbative nature.  Moreover, the semi-local version of the  Almost reducibility conjecture has no reasonable equivalent in the smooth (or even Gevrey) setting \cite{AK3}.  Still, one can ask:
\begin{question} Is the semi-local  version of the Almost reducibility conjecture true  for cocyles in  quasi-analytic classes?
\end{question}

Let's say that a cocycle is {\it stable} if it is not accumulated by non-uniformly hyperbolic systems (with the same frequency vector on the base). Having in mind Avila's classification one can ask:
\begin{question}Is every stable  cocycle  in $SW_{}^k(\T^d,SL(2,\R))$, $k=\infty,\omega$,  almost-reducible? 
\end{question}
\subsection{The symplectic case}
Cocycles in $SW^k_{}(\T^d,Sp(2n,\R))$ are of interst when one tries to understand the dynamics of a symplectic diffeomorphism in the neighborhood of an invariant torus (they appear as linearized dynamics) or in the study of quasi-periodic Schr\"odinger operators on strips $\Z\times\{1,\ldots,n\}$. For such cocycles one can define $2n$ Lyapunov exponents (symmetric with respect to 0) and one fibered Maslov index which plays the role of a fibered rotation number (cf. \cite{Xu1} and the references there).

We denote by $SO(2,\R)$ the set of symplectic rotations $R_{t}=\bm (\cos t)I_{n} &-(\sin t) I_{n}\\ (\sin t) I_{n}&(\cos t)I_{n}\em$.
\begin{question} Let $(\a,A)\in SW^\infty(\T,Sp(2n,\R))$ homotopic (resp. non homotopic) to the identity   where $\a\in\T$ is   (recurrent) diophantine.
Is it true that for Lebesgue almost all $t\in\R$ the following dichotomy holds: either the cocycle $(\a,R_{t}A)$ is $C^\infty$-reducible (resp. $SO(2,\R)$-reducible) or its upper Lyapunov exponent is positive? 
\end{question}
When $n=1$ the answer is positive (\cite{AK-annals} for the case homotopic to the identity, \cite{AK-inventiones} for the case non-homotopic to the identity). The proof of this result is based on a renormalization procedure which works when the cocycle has some mild boundedness property  and on a reduction to this case based on Kotani theory.  In the case $n\geq 2$ such a Kotani theory was developed by  Xu in  \cite{Xu1}, \cite{Xu2}. Following the same strategy as in \cite{AK-annals} one should be  then reduced to studying cocycles with values in the maximal compact subgroup of $Sp(2n,\R)$. Unfortunately, one cannot conclude like in the case $n=1$ since no reasonable {\it a priori} notion of fibered rotation number can be defined  for cocycles with values in non-abelian compact groups (they can be defined {\it a posteriori} once one knows the cocycle is reducible; see  \cite{Karaliolios-JMD}, \cite{Karaliolios-SMF} for related  results).
\subsection{The case $G={\rm Diff}_{0}^\infty(\T)$}
A cocycle $(\a,A)\in SW(\T^d,SL(2,\R))$ naturally produces a {\it projective cocycle} $(\a,\bar A)\in SW(\T^d,{\rm Hom}(\bS^1))$ where ${\rm Hom}(\bS^1)$ is the group of homographies acting on $\bS^1$; namely $\bar A(x)\cdot v=(A(x)v)/\|A(x)v\|$. It is thus natural to look at the more general case where the underlying group is the group of orientation preserving diffeomorphisms of the circle. In that case one can still define a fibered rotation number \cite{Herman83}. For the  topological aspects of the  theory of such quasi-periodically forced circle diffeomorphisms see \cite{BeJa}.
\begin{question}[Non-linear Eliasson Theorem] Let $\a\in\T^d$ be a fixed diophantine vector and $G={\rm Diff}_{0}^\infty(\R/\Z)$. Does there exist $k_{0},\e_{0}$ depending only on $\a$ such that for any $(\a,A)\in SW^\infty(\T^d,{\rm Diff}_{0}^\infty(\R/\Z))$ of the form $(\a,A)(x,y)=(x+\a,y+\beta+f(x,y))$ with $\|f\|_{C^{k_{0}}}\leq \e_{0}$ and $\rho(\a,A)$ diophantine, the cocycle $(\a,A)$ is $C^\infty$-reducible?
\end{question}
 When $G={\rm Hom}(\bS^1)$ the answer is positive and  is (the $C^\infty$-version of) a theorem of Eliasson \cite{El-CMP92} which has  many consequences in the theory of quasi-periodic Schr\"odinger operators.  
If one allows $\e_{0}$ to depend on  $\rho$ then the result is true and is essentially a (generalization of a)  theorem by Arnold. Its proof is classical KAM theory.  In \cite{KWYZ} a result of rotations-reducibility is proved where $\e_{0}$ depends  on $\rho$ but with considerably weaker assumption on $\a$ than KAM theory usually allows (compare with \cite{AFK}, \cite{rb6} for stronger  results in the case of  linear cocycles).

\section{Mixing surface flows.}

\subsection{Spectral type.} Area preserving surface flows provide the lowest dimensional setting in which it is interesting to study conservative systems. 
Such flows are sometimes called multi-valued Hamiltonian flows to emphasize their relation with solid state physics that was pointed out by  Novikov \cite{novikov}.  
Via Poincaré sections, these flows are related to special flows above circle rotations or more generally above IETs (Interval exchange transformations). One can thus view them as time changes of translation flows on surfaces.  

Katok and then Kochergin showed the absence of mixing of area preserving flows on the two torus if they do not have singularities \cite{katok,Koc1}.

The simplest mixing examples are those with one (degenerate) singularity on the two torus produced by Kochergin in the 1970s \cite{Koc2}.  Kochergin flows are time changes of linear flows on the two torus with an irrational slope and with a rest point (see Figure 1).

Multi-valued Hamiltonian flows on higher genus surfaces can also be mixing (or mixing on an open ergodic component) in the presence of non-degenerate saddle type singularities that have some asymmetry  (see Figure 2).
Such flows are called Arnol'd flows and their mixing property, conjectured by Arnol'd in \cite{arnold},  was obtained by Khanin-Siani \cite{SK} and in more generality by Kochergin \cite{Kc3}.  
 Note that Ulcigrai proved in \cite{CU3} that  area preserving flows with non-degenerate saddle singularities are generically not mixing (due to symmetry in the saddles).

\begin{figure}
\begin{center} 
\begin{tabular}{cc}
   \resizebox{!}{3cm}{ \includegraphics{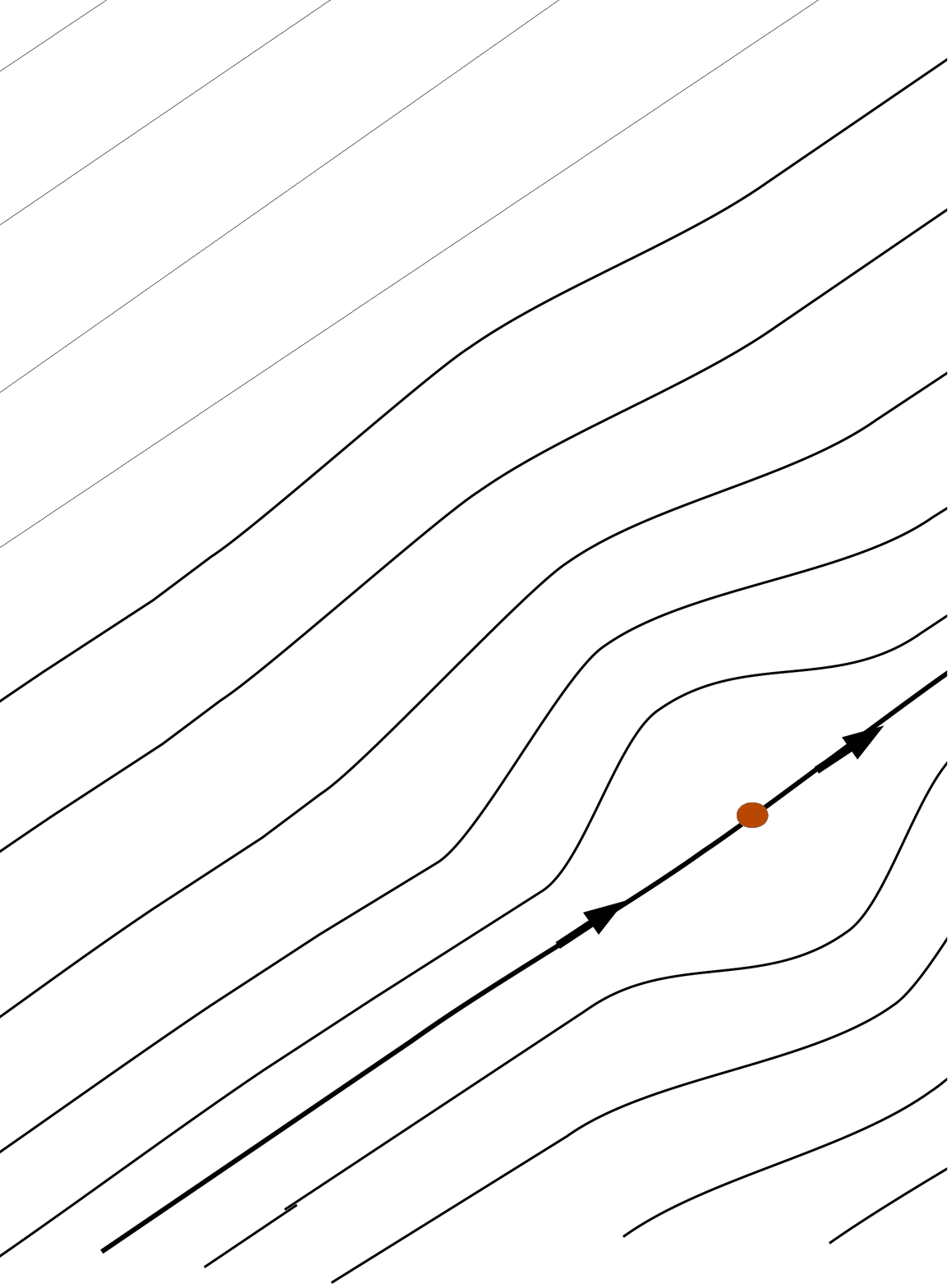}}   &   \resizebox{!}{3cm}{ \includegraphics{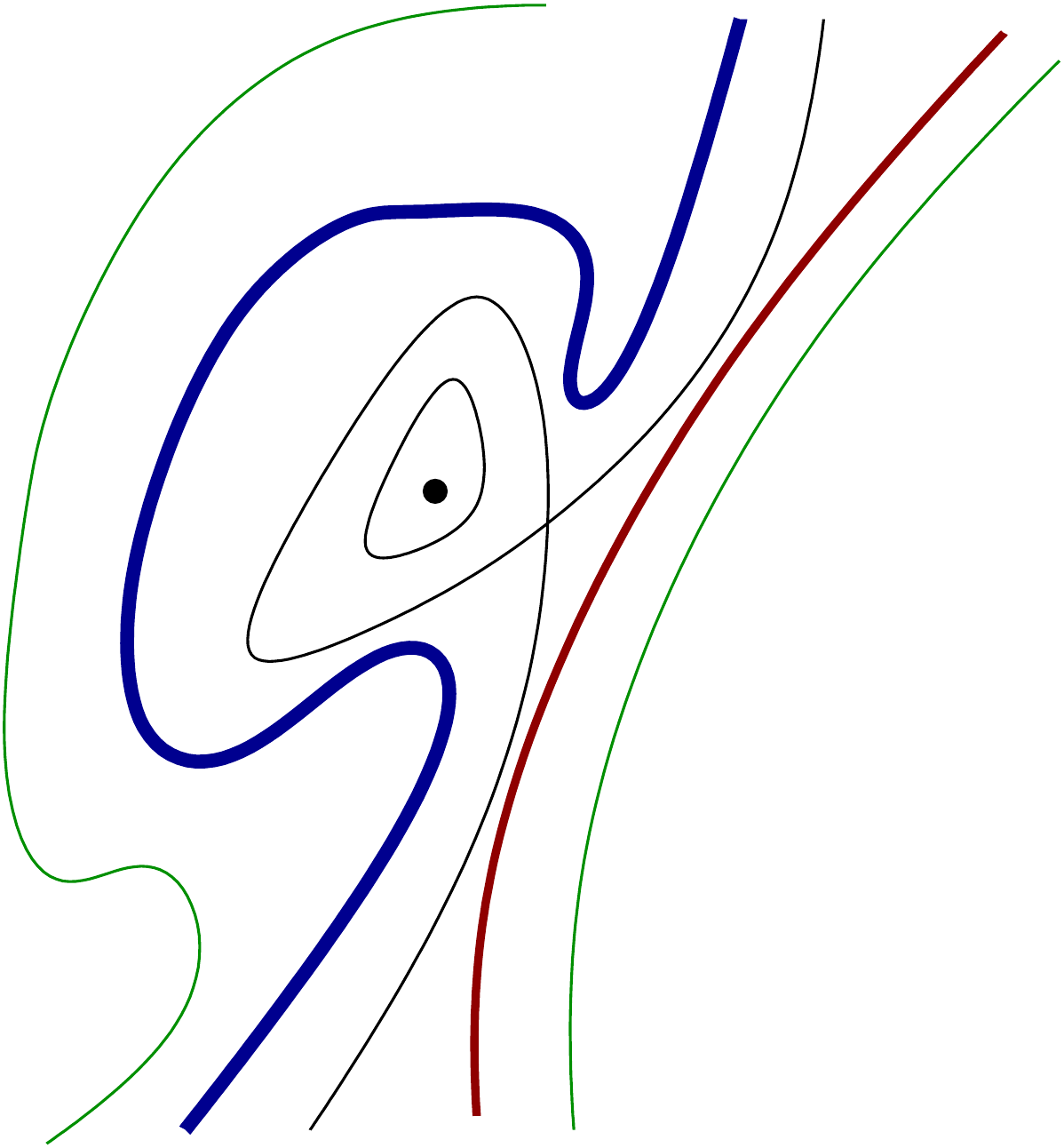} }
     \\[\abovecaptionskip]    \resizebox{!}{0.15cm}{Fig 1.  \small Degenerate saddle acting as a stopping point}  &     \resizebox{!}{0.15cm}{Fig 2. \small   Non-degenerate saddle that causes asymmetry}   \\
\end{tabular}
\end{center} 
       
\end{figure}

\begin{question} Study the spectral type and spectral multiplicity of mixing flows on surfaces. 
\end{question}
By spectral type of a flow $\{T^t\}$ we mean the spectral type of the associated Koopman operator $U_t : L^2(M,\mu)  : f \to f \circ T^t$.

\begin{figure}
\begin{center} 
\end{center} 
       
\end{figure}

It was proved in \cite{FFK} that Kochergin flows with a sufficiently strong power like singularity have for almost every slope a maximal spectral type that is equivalent to Lebesgue measure. 
The study of the spectral multiplicity of these flows is interesting in its relation to the Banach problem on the existence of a dynamical system with simple Lebesgue spectrum. It is probable however that the spectral multiplicity of Kochergin flows is infinite. Mixing reparametrizations of linear flows with simple spectrum were obtained in \cite{rankone} and it would be interesting to study their maximal spectral type following \cite{FFK}.

\begin{question} Is it true that Arnol'd mixing flows have in general a purely singular spectral type?
\end{question}

Arnol'd conjectured a power-like decay of correlation in the non-degenerate asymmetric case, but the decay is more likely to be logarithmic, at least between general regular observables or characteristic functions of regular sets such as balls or squares. Even a lower bound on the decay of correlations is not sufficient to preclude absolute continuity of the maximal spectral type. However, an approach based on slowly coalescent periodic approximations as in \cite{Fa3} may be explored in the aim of proving that the spectrum is purely singular.

\subsection{Spectral type of related systems}

\begin{question} Prove that all IET have a purely singular maximal spectral type.
\end{question}
It is known that almost every IET, namely those that are not of constant type, are rigid. It follows that their maximal spectral type is purely singular. For the remaining IETs, partial rigidity was proven by Katok and used to show the absence of mixing, but proving that the spectral type is purely singular appears to be more delicate.

\begin{question} Prove that on $\T^3$ there exists a real analytic strictly positive reparametrization of a minimal translation flow that has a Lebesgue maximal spectral type. 
\end{question}

The difference with the Kochergin flows is that such flows would also be uniquely ergodic.
Mixing real analytic reparametrizations of linear flows on $\T^3$ were obtained in \cite{etds2}.

\subsection{Multiple mixing.} 

The question of multiple mixing for mixing systems is one of the oldest unsolved questions of ergodic theory.
\begin{question} Are all mixing surface flows mixing of all orders? 
\end{question}

Arnold and Kochergin mixing conservative flows on surfaces stand
as the main and almost only natural class of mixing transformations for which
higher order mixing has not been established nor disproved in full generality. Under suitable
arithmetic conditions on their unique rotation vector, of full Lebesgue measure
in the first case and of full Hausdorff dimension in the second, it was shown in \cite{FKa}  that
these flows are mixing of any order,  \cite{KKU} for flows on higher genus surfaces).

\section{Ergodic theory of diagonal actions on the space of lattices and applications to metric Diophantine approximation} 

\medskip
 
 The Diophantine properties of linear forms of one or several variables evaluated at integer points are intimately related to the divergence rates of some orbits under some diagonal actions in the space of (linear or affine) lattices of $\R^n$. This link is due to what can be called the Dani correspondence principle between the small values of the linear forms  on one hand and the visits to the cusp of certain orbits of certain diagonal actions on the space of lattices (affine lattices in the case of inhomogeneous linear forms). The ergodic study of diagonal and unipotent actions on the space of lattices provides indeed an efficient substitute to the continued fraction algorithm that played a crucial role in the rich metric theory of Diophantine approximations in dimension 1. There is a number of important contributions to number theory related to this principle and to progress in the theory of homogeneous actions for example the surveys  \cite{dani,kleinbock,EL,eskin,M-LH,M2}).    We mention here a list of questions related to the statistical properties of Kronecker sequences that can be approached using this same principle. More details and questions can be found in \cite{torusrev}.

\subsection{Kronecker sequences} A quantitative measure of uniform distribution of Kronecker sequences is given by the discrepancy function:  for a set 
$\cC \subset \T^d$ let 
$$D(\alpha,x,\cC, N)=\sum_{n=0}^{N-1} \one_{\cC}(x+k\alpha)-N {\rm volume}(\cC)$$
where $(\a,x) \in \T^d \times \T^d$ and  $\one_\cC$ is the 
characteristic function of the set $\cC$.

Uniform distribution of the sequence $x+k\a $ on $\T^d$ is equivalent to the fact that, 
for regular sets $\cC,$ 
$D(\alpha,x,\cC, N)/N \to 0$ as $N \to \infty$. A step further  is the study of the 
rate of convergence to $0$ of $D(\alpha,x,\cC, N)/N$.

Already for $d=1$, it is clear that if $\a \in \T-\Q$ is fixed, the discrepancy $D(\alpha,x,\cC, N)$ displays an oscillatory 
behavior according to the position of $N$ with respect to the denominators of the best rational approximations of  $\a$. 
A great deal of work in Diophantine approximation has been done on giving upper and lower bounds to the oscillations of the discrepancy function (as a function of $N$) in relation with 
the arithmetic properties of   $\a \in \T^d$.

 In particular, let
\begin{equation*}
\overline{D}(\a, N)=\sup_{\Omega \in \mathbb{B}} D(\alpha,0,\Omega,N)
\end{equation*}
where the supremum is taken over all sets $\Omega$ in some natural class of sets $\mathbb{B}$, for example balls or boxes.

The case of (straight) boxes was extensively studied, and properties of the sequence $\overline{D}(\a, N)$ were obtained  with a special emphasis on their relations with the Diophantine approximation properties of $\a.$ 
In particular, 
 \cite{beck1} proves that when $\mathbb{B}$ is the set of straight boxes in $\T^d$ then for arbitrary positive increasing function $\phi(n)$ 
\begin{equation} \sum_{n} \frac{1}{\phi(n)}<\infty  \iff \frac{\overline{D}(\a, N)}{(\ln N)^d \phi(\ln\ln N)}   \begin{array}{l}
\text{ is  bounded for}\\
\text{ almost every  }\a \in \T^d.
\end{array} 
   \label{beckbound}  \end{equation}

In dimension $d=1$, this result 
is the content of Khinchine theorems obtained in the early 1920's,
and it follows easily from well-known  results from 
the metrical theory of continued fractions (see for example the introduction of \cite{beck1}).
The higher dimensional case is significantly more difficult and many questions that are relatively easy to settle in dimension $1$ remain open. We mention some here and refer to \cite{beck1,Nie} for others. 
 
\begin{question}
Is it true that $\limsup\frac{\overline{D}(\alpha,N)}{\ln^d N}>0$ {\bf for all} $\alpha\in \T^d$?
\end{question}

\begin{question}
Is it true that there exists $\a$ such that $\limsup\frac{\overline{D}(\alpha,N)}{\ln^d N}<+\infty$?
\end{question}

The above questions and results can be asked for balls and more general convex sets. 
\begin{question} Is it true that for any $\eps>0$, for almost every $\a \in \T^d$ and for any convex set $\cC$ in $\T^d$
$$\frac{D(\alpha,0,\cC, N)}{ N^{\frac{d-1}{2d}+\eps}}$$ is bounded? 
\end{question}

 The bound in (\ref{beckbound}) 
 focuses on how bad can the discrepancy become along a subsequence of $N$, for a fixed $\a$ in a full measure set.  In a sense, it  deals with the worst case scenario and do not capture the oscillations of the discrepancy.

Another point of view is to let $(\a,x)\in \T^d \times \T^d$ be random and have limit laws that hold for {\it all} $N$. By random we mean distributed according to a smooth density on the tori. For $d=1$, this was done by 
Kesten who proved in the 1960s that the discrepancies of the number of visits of the Kronecker sequence to an interval, normalized by $\rho \ln N$ (where $\rho$ depends on the interval but is constant if the length of the interval is irrational) converges to a Cauchy distribution. 

One can ask whether Kesten's convergence remains valid for a fixed $x$. Another question is what happens in higher dimension? In particular : 

\begin{question} Is it true that there exists $\rho>0$, such that when $\cC$ is a generic box in $\T^d$ and $\a$ is uniformly distributed on $\T^d$,  then $\frac{D(\alpha,0,\cC, N)}{\rho(\ln N)^d}$  converges in distribution to the Cauchy law? 
\end{question}

In \cite{db_poisson} this was proved when $x$ and the box $\cC$ are also random (a shape is randomized by applying small deformations distributed according to a smooth measure on the space of isometries).  It was shown in \cite{db_convex} that in the case of a strictly convex shape $\cC \subset \T^d$ one has $ \frac{D(\a,x,r\cC,N)}{r^{\frac{d-1}{2}} N^{\frac{d-1}{2d}} }$ converges in distribution to a non standard law when $(\a,x) \in \T^d \times \T^d$ and $r>0$ are random. The convex set $r\cC$ is the rescaled set from $\cC$ by factor $r$ around some fixed point inside $\cC$.

A {\it semialgebraic} set $\cC$ in $\T^d$ is a set defined by a finite number of algebraic inequalities. This includes a diverse collection of
sets such as balls, cubes, cylinders, simplexes etc. Following \cite{db_convex} we ask 

\begin{question} 
Assume $\cC$ is semialgebraic. Does there exist a sequence $a_N=a_N(\cC)$ such that 
 $ \frac{D(\a,x,\cC,N)}{a_N}$ converges in distribution when $(x, \alpha) \in \times \T^d \times \T^d$ are random.
\end{question}

One can study the fluctuations of the ergodic sums above toral translations for functions other from characteristic functions. The following is interesting for its connection with number theory as well as with the ergodic theory of some natural classes of dynamical systems such as surface flows.
 
\begin{question}  Study  the behavior of the ergodic sums $\sum_{n_1=1}^N  A(x+n \alpha)$ for functions $A$ that are smooth except for a finite number of singularities.  
\end{question}

The fluctuations can be studied for fixed $\a$ or $x$, as well as for random values. One should then try to classify the fluctuations according to the type of the singularities : power, fractional power, logarithmic (we refer to \cite{M2,torusrev} for more details and questions).

\subsection{Higher dimensional actions.} Replacing the $\Z$ action by translation with $\Z^k$ actions (see we get  following \cite{torusrev}

\begin{question} 
Study the ergodic sums $\sum_{j=1}^m \sum_{n_j=1}^N A(x+\sum_{j=1}^m  \alpha_j n_j)$, with $(x,\a_1,\ldots,\a_m) \in (\T^d)^{m+1}$. 
\end{question}

In the case where $A=\chi_I-|I|$ and $\chi_I$  the indicator of an interval we get the following possible extension of Kesten's theorem to the statistical behavior of linear forms.
\begin{question} Show that as $x \in \T$ and $\a \in \T^m$ are random $$\frac{1}{\rho(\ln N)^d}\sum_{j=1}^m \sum_{n_j=1}^N A(x+\sum_{j=1}^m  \alpha_j n_j)$$
 converges in distribution to a Cauchy law for some $\rho>0$.  \label{linearform}
\end{question}

One can also investigate analogues of the Shrinking Targets Theorems of \cite{DFV} for $\Z^k$ actions. 
\begin{question}
Let $l, \hat l: \R^d\to \R, $ be linear forms with random coefficients, $Q:\R^d\to \R$ be a positive definite quadratic form.
Investigate limit theorems, after adequate renormalization, for the number of solutions to

(a) $\{l(n)\}Q(n) \leq c, |n|\leq N;$ 

(b) $\{l(n)\} |\hat l(n)|\leq c, |n|\leq N;$

(c) $|l(n)Q(n)|\leq c, |n|\leq N;$ 

(d) $|l(n)\hat l(n)|<c, |n|\leq N.$
\end{question}

\end{document}